\baselineskip=16pt

\font\twelvebf=cmbx12

\bigskip
\bigskip

\def\C{{\bf C}}
\def\Ch{{\C[[h]]}}
\def\N{{\bf N}}

\def\e{\eqno}
\def\l{\ldots}
\def\n{\noindent}
\def\Uq{U_q[sl(n+1)]}

\def\a{\tilde a}
\def\b{\tilde b}

\def\q{{\bar q}}
\def\s{\smallskip}
\def\b{\bigskip}

\def\>{\rangle}
\def\<{\langle}
\def\v{\varphi}

\hfill {{\twelvebf Preprint UR-1522 $\;\;$ 3/98}

\vskip 4cm

\noindent
{\twelvebf A Holstein-Primakoff and a Dyson realization for the

\noindent
quantum algebra U$_q$[sl(n+1)]  }

\leftskip 50pt
\vskip 32pt
\noindent
{\bf Tchavdar D. Palev}\footnote*{Permanent address: 
Institute for Nuclear Research and Nuclear Energy, 1784 Sofia, 
Bulgaria; $e$-mail: tpalev@inrne.acad.bg}

\noindent
Department of Physics and Astronomy, University of Rochester,
Rochester, New York 14627
\vskip 12pt

%\vfill\eject 

\noindent 
{\bf Abstract.} The known Holstein-Primakoff and Dyson realizations of
the Lie algebra $sl(n+1),\;n=1,2,\ldots$ in terms of Bose operators
(Okubo S 1975 {\it J. Math. Phys.} {\bf 16} 528) are generalized to
the class of the quantum algebras  $U_q[sl(n+1)]$ for any $n$.
It is shown how the elements of $U_q[sl(n+1)]$ can be expressed
via $n$ pairs of Bose creation and annihilation operators.

%\vfill\eject
\leftskip 0pt
\vskip 48pt

\vskip 12pt

\bigskip
\noindent
In the present note we write down an analogue of the Dyson (D) 
realization and of the Holstein-Primakoff (HP)
realization  for the quantum algebra
$\Uq$. Initially the HP and the D realizations were given for $sl(2)$
[1, 2]. The generalization for $gl(n+1)$ is due to Okubo [3].
In the realization of Okubo the elements of $gl(n+1)$ are
expressed as functions of $n$ pairs of Bose creation and annihilation
operators (CAOs), namely operators $a_1^\pm,a_2^\pm,\l,a_n^\pm$, which
satisfy the known commutation relations
$$
[a_i^-,a_j^+]=\delta_{ij},\quad [a_i^+,a_j^+]=[a_i^-,a_j^-]=0,
\quad i,j=1,\ldots,n.  \eqno(1)
$$
This realization is ``more economical'' than the known 
Jordan-Schwinger realization, which expresses $gl(n+1)$ 
via $n+1$ pairs of Bose CAOs.

The motivation in the present work stems from the various applications
of the Holstein-Primakoff and of the Dyson realizations in theoretical
physics. Beginning with [1] and [2] the HP and D realizations were
constantly used in condensed matter physics. Some other early
applications can be found in the book of Kittel [4] (more recent
results are contained in [5]). For applications in nuclear physics see
[6, 7] and the references therein, but there are, certainly, several
other publications. In view of the importance of the quantum algebras
for physics, one could expect that the generalization of the Dyson and
of the Holstein-Primakoff realizations to the the case of quantum
algebras may be of interest too. This is in fact the case for the only
known so far $q-$analogs of the HP realization, namely those of
$U_q[sl(2)]$ and $U_q[sl(3)]$ [8-15].

To begin with we recall the definition of $U_q[sl(n+1)]$
in the sense of Drinfeld [16]. Let $\Ch$ be
the complex algebra of the formal power series in the indeterminate 
$h$,
$q=e^{h/2}\in \Ch$. Then $U_q[sl(n+1)]$ is a Hopf algebra, which is a
topologically free $\Ch$ module (complete in the $h-$adic topology),
with generators $\{h_i,e_i,f_i\}_{i=1,\l,n}$ and

\s\n
1. Cartan relations
$$
\eqalignno{
&  [h_i,e_j]=(2\delta_{ij}-\delta_{i,j-1}-\delta_{i-1,j})e_j,& (2a) 
\cr
%&&\cr
&  [h_i,f_j]=-(2\delta_{ij}-\delta_{i,j-1}-\delta_{i-1,j})f_j,& 
(2b)\cr   
%&&\cr
&  [e_i,f_j]=\delta_{ij}{{q^{h_i}-\q^{h_i}}\over{q-{\bar q}}}.& (2c) 
\cr
}
$$

\n
2. Serre relations
$$
\eqalignno{
& [e_i,e_j]=0, \quad  [f_i,f_j]=0,\quad \vert i-j \vert \neq 1,& 
(3a)\cr
%&&\cr
& [e_i,[e_i,e_{i \pm 1}]_{\bar q}]_q=
[e_i,[e_i,e_{i \pm 1}]_q]_{\bar q}=0, &(3b) \cr
%&&\cr
& [f_i,[f_i,f_{i \pm 1}]_{\bar q}]_q=
[f_i,[f_i,f_{i \pm 1}]_q]_{\bar q}=0. & (3c) \cr 
}
$$
Above and throughout $[a,b]=ab-ba$,  $[a,b]_x=ab-xba$, $\q=q^{-1}$.
We do not write the other Hopf algebra maps $(\Delta,\; \varepsilon,
S)$, since we will not use them. They are certainly also a part of the
definition.

The Dyson and the Holstein-Primakoff realizations are different
embeddings of $\Uq$ into the Weyl algebra $W(n)$. We define the latter
as a topologically free $\Ch$ module and an associative unital algebra
with generators $a_1^\pm,\l,a_n^\pm$ and relations (1).

\s\n
{\it Remark.} In the physical applications it is often more convenient
to consider $h$ and $q$ as complex numbers, $h, q\in \C$. Then all our
considerations remain true provided $q$ is not a root of 1. The
replacement of $q\in \Ch$ with a number corresponds to a factorization
of $\Uq$ and $W(n)$ with respect to the ideals generated by the
relation $q=number$. The factor-algebras $U_q[sl(n+1)]$ and $W(n)$ are
complex associative algebras. However the completion in the $h$-adic
topology has left a relevant trace: after the factorization the
elements of $U_q[sl(n+1)]$ and of $W(n)$ are not only polynomials of
their generators. In particular the functions of the CAOs, which
appear in the D and in the HP realizations (see (4) and (10) bellow)
are well defined as elements from $W(n)$.

Now we are ready to state our main results. Let
$[x]={{q^x-\q^x}\over{q-\q}}$, 
$N_i=a_i^+a_i^-$ and $N=N_1+\l+N_n$.

\smallskip\noindent
{\it Proposition 1 (Dyson realization).} The linear map $\varphi:
\Uq \rightarrow W(n)$, defined on the generators as
$$
\eqalignno{
&
\v(h_1)=p-N-N_1,\quad \v(h_i)=N_{i-1}-N_i,\quad i=2,3,\ldots,n,
& (4a)   \cr
&&\cr
&
\v(e_1)={[N_1+1]\over{N_1+1}}[p-N]b_1^-,\quad
\v(e_i)={[N_i+1]\over{N_i+1}}b_i^-b_{i-1}^+,
\quad i=2,\ldots,n, & (4b) \cr
&&\cr
&
\v(f_1)=b_1^+,\quad
\v(f_i)={[N_{i-1}+1]\over{N_{i-1}+1}}b_i^+b_{i-1}^-,
\quad i=2,\ldots,n, & (4c) \cr
}
$$
is a morphism of $\Uq$ into $W(n)$ for any $p \in \C$.

\s
The proof is straightforward.
In the intermediate computations the following 
relation is useful:
$$
f(N_1,\l,N_i,\l,N_n)a_j^\pm = 
a_j^\pm 
f(N_1\pm\delta_{1j},\l,N_i\pm\delta_{ij},\l,N_n\pm\delta_{nj}),\e(5)
$$
where $f(N_1,\l,N_i,\l,N_n) \in W(n)$  is a function of the
number operators $N_1,\l,N_i,\l,N_n$.

We have derived the D realization (4) on the ground of an alternative
to the Chevalley definition of $\Uq$ [17]. This derivation
together with the expressions for (the analogs of) all Cartan-Weyl
generators via CAOs will be given elsewhere.

Similar as for $sl(n+1)$, the D realization defines 
an infinite-dimensional representation of $\Uq$  in the Fock space
${\cal F}(n)$ with an orthonormed basis
$$
|l\>\equiv |l_1,\ldots,l_n\>={(a_1^+)^{l_1}\ldots (a_n^+)^{l_n}\over
\sqrt{l_1!\ldots l_n!}}|0\>,\quad l_1,\ldots, l_n=0,1,2,\l. \eqno(6)
$$
If $p$ is a positive integer, $p\in \N$, the representation is
indecomposible: the subspace 
$$
{\cal F}_1(p;n)=lin.env.\{|l_1,\ldots,l_n\>|l_1+\l+l_n> p\}\e(7)
$$
is an invariant subspace, whereas its orthogonal
compliment
$$
{\cal F}_0(p;n)=lin.env.\{|l_1,\ldots,l_n\>|l_1+\l+l_n\le p\}\e(8)
$$
is not an invariant subspace. If $p\notin \N$, the representation is
irreducible. In all cases however, and this is the disadvantage of the
D realization, the representation of $\Uq$ in ${\cal F}(n)$ is not
unitarizable with respect to the antilinear anti-involution
$\omega:\Uq \rightarrow \Uq$, defined on the generators as

$$
\omega(h_i)=h_i,\quad \omega(e_i)=f_i,\quad i=1,\l,n.\e(9)
$$

In order to ``cure'' this disadvantage we pass to introduce the HP
realization.

\smallskip\noindent
{\it Proposition 2 (Holstein-Primakoff realization).} The linear map
$\pi:
\Uq \rightarrow W(n)$, defined on the generators as
$$
\eqalignno{
& 
\pi(h_1)=p-N-N_1,\quad \pi(h_i)=N_{i-1}-N_i,\quad i=2,3,\ldots,n,
& (10a)   \cr
&&\cr
& 
\pi(e_1)=\sqrt{{[N_1+1]\over{N_1+1}}[p-N]}\;a_1^-,\quad
\pi(e_i)=
\sqrt{{[N_{i-1}]\over{N_{i-1}}}{[N_i+1]\over{N_i+1}}}\;a_i^-a_{i-1}^+,
\quad i=2,3,\ldots,n,  & (10b) \cr
&&\cr
& 
\pi(f_1)=\sqrt{{[N_1]\over{N_1}}[p-N+1]}\;a_1^+,\quad
\pi(f_i)=\sqrt{{[N_{i-1}+1]\over{N_{i-
1}+1}}{[N_i]\over{N_i}}}\;a_i^+a_{i-1}^-,
\quad i=2,3,\ldots,n, & (10c) \cr
}
$$

\s\n
is a morphism of $\Uq$ into $W(n)$ for any $p \in \C$. If $p\in\N$,
then ${\cal F}_0(p;n)$ and ${\cal F}_1(p;n)$ are invariant
subspaces; ${\cal F}_0(p;n)$ carries a finite-dimensional
irreducible representation; it is unitarizable with
respect to the anti-involution (9) and the metric defined with the
orthonormed basis (6), provided $q>0$.

The proof is straightforward: the verification of the defining
relations (2) and (3) can be carried out on a purely algebraic
level. The circumstance that  ${\cal F}(n)$ is a direct sum of
its invariant subspaces ${\cal F}_0(p;n)$ and ${\cal F}_1(p;n)$
is due to the the factor $\sqrt{[p-N]}$ in (10b) and
$\sqrt{[p-N+1]}$ in (10c). 
If $q>0$, then ($(\;,\;)$ denotes the scalar product)
$$
(\pi(h_i)|l\>,|l'\>)=(|l\>,\pi(h_i)|l'\>),\quad
(\pi(e_i)|l\>,|l'\>)=(|l\>,\pi(f_i)|l'\>) \quad
\forall \;|l\>, |l'\>\in {\cal F}_0(p;n),\;\;i=1,\l,n.
$$
Therefore the representation of  $\Uq$ in
${\cal F}_0(p;n)$ is unitarizable.

Let us note that the HP realization (10) of $\Uq$ can be
easily expressed also in terms of deformed oscillator operators
$\a_i^\pm,\;{\tilde N}_i$, $i=1,\l,n$, namely operators which 
satisfy the relations [19-21]:
$$
[\a_i^-,\a_j^+]_q=\delta_{ij} q^{-{\tilde N}_i},\quad
[{\tilde N}_i,\a_j^\pm]=\pm \delta_{ij}\a_j^\pm,\quad
[\a_i^\pm,\a_k^\pm]=[{\tilde N}_i,{\tilde N}_k]=0,\quad
i\ne k. \e(11)
$$ 
From (10) and the relations between the deformed and the nondeformed
operators [22]
$$
\a_i^-=\sqrt{{[N_i+1]\over{N_i+1}}}\;a_i^-,\quad
\a_i^+=\sqrt{{[N_i]\over{N_i}}}\;a_i^-,\quad
{\tilde N}_i=N_i,\e(12)
$$
one obtains a $q-$analogue of the HP realization for any $n$:
$$
\eqalignno{
& 
\pi(h_1)=p-{\tilde N}-{\tilde N}_1,
\quad \pi(h_i)={\tilde N}_{i-1}-{\tilde N}_i,\quad i=2,3,\ldots,n,
& (13a)   \cr
%&&\cr
& 
\pi(e_1)=\sqrt{[p-{\tilde N}]}\;\a_1^-,\quad
\pi(e_i)=\a_i^-\a_{i-1}^+,
\quad i=2,3,\ldots,n,  & (13b) \cr
%&&\cr
& 
\pi(f_1)=\sqrt{[p-{\tilde N}+1]}\;\a_1^+,\quad
\pi(f_i)=\a_i^+\a_{i-1}^-,
\quad i=2,3,\ldots,n, & (13c) \cr
}
$$
To the best of our knowledge such $q-$deformed analogs of HP
realizations are available so far only for $U_q[sl(2)]$ [8-14] and for
$U_q[sl(3)]$ [15].

Finally, adding to the generators of $\Uq$ an additional central
element $I$, and setting \break 
$\varphi(I)=\pi(I)=p$, one obtains a D
and a HP realizations of $U_q[gl(n+1)]$.

\b\b\n
{\bf Acknowledgments }

\b\n
The author is thankful to Prof. S. Okubo for the kind invitation to
conduct the research under the Fulbright Program in the Department of
Physics and Astronomy, University of Rochester. This work was
supported by the Fulbright Program, Grant No 21857.

\vfill\eject

\noindent
{\bf References}

\vskip 12pt

\+ [1] & Holstein T and Primakoff H 1949 
         {\it Phys. Rev.} {\bf 58} 1098 \cr            

\+ [2] & Dyson F J 1956 {\it Phys. Rev.} {\bf 102} 1217 \cr

\+ [3] & Okubo S 1975 {\it J. Math. Phys.} {\bf 16} 528 \cr

\+ [4] & Kittel C 1063 {\it Quantum Theory of Solids} (Willey, 
          New York) \cr

\+ [5] & Caspers W J 1989  {\it Spin Systems} (World Sci. Pub. Co., 
Inc.,
          Teanek) \cr

\+ [6] & Klein A and Marshalek E R 1991
         {\it Rev. Mod. Phys.} {\bf 63} 375 \cr

\+ [7] & Ring P and Schuck P  {\it The Nuclear Mani-Body Problem}
          (Springer-Verlag, New York, Heidelberg, Berlin) \cr

\+ [8] & Chaichian M, Ellinas D and Kulish P P 1990 
         {\it Phys. Rev. Lett.} {\bf 65} 980  \cr

\+ [9] & Quesne C  1991 {\it Phys. Lett. A} {\bf 153}
         303 \cr

\+ [10] &$\;$ Chakrabarti R and Jagannathan R 1991 
         {\it J. Phys.  A~: Math.\ Gen.} 
         {\bf 24}  L711 \cr

\+ [11] &$\;$ Katriel J and Solomon A I 1991 
         {\it J. Phys. \ A~: Math.\ Gen.} {\bf 24} 2093 \cr

\+ [12] &$\;$ Yu Z R 1991 
         {\it J. Phys. \ A~: Math.\ Gen.} {\bf 24} L1321 \cr

\+ [13] &$\;$ Kundu A and Basu Mallich B 1991 {\it Phys. Lett. A} 
          {\bf 156} 175  \cr

\+ [14] &$\;$ Pan F 1991 {\it Chin. Phys. Lett.} {\bf 8} 56 \cr

\+ [15] &$\;$ da-Providencia J 1993 
         {\it J. Phys. \ A~: Math.\ Gen.} {\bf 26} 5845 \cr

\+ [16] &$\;$ Drinfeld V 1986 {\it Quantum Groupd (Proc. Int. 
         Congress of Mathematics (Berkeley, 1986)) }\cr
\+      & ed A M Gleasom (Providence, RI: American Physical
          Society) p 798  \cr

\+ [17] &$\;$ Palev T D and Parashar P 1998 {\it Lett.Math. Phys.} 
         {\bf 43} 7 \cr

\+ [18] &$\;$ Polychronakos A P 1990 {\it Mod.Phys. Lett.}
          {\bf A5} 2325 \cr  

\+ [19] &$\;$ Macfarlane A J 1989 {\it J.\ Phys.\ A~: Math.\ Gen.}
          {\bf 22}  4581  \cr

\+ [20] & $\;$ Biedenharn L C 1989 {\it J.\ Phys.\ A~: Math.\ Gen.} 
          {\bf 22}  L873 \cr

\+ [21] & $\;$ Sun C P and Fu H C 1989 {\it J.\ Phys.\ A~: Math.\ 
Gen.} 
         {\bf 22}  L983  \cr

\+ [22] &$\;$ Polychronakos A P 1990 {\it Mod.Phys. Lett.}
            {\bf 5} 2325\cr

\end